\documentclass[12pt]{article}
\usepackage{wasysym,amssymb,eufrak,indentfirst,cite,bibspacing,color}
\begin{document}
\baselineskip 17pt
\title{\textbf{On partial $\Pi$-property of subgroups of finite groups}\thanks{Research is supported by a NNSF grant of China (grant \#11071229) and Research Fund for the Doctoral Program of Higher Education of China (Grant 20113402110036).}}
\author{Xiaoyu Chen, Wenbin Guo\\
{\small Department of Mathematics, University of Science and Technology of China,}\\ {\small Hefei 230026, P. R. China}\\
 {\small E-mail: jelly@mail.ustc.edu.cn, $\,$wbguo@ustc.edu.cn}
}
\date{}
\maketitle
\begin{abstract}
Let $H$ be a subgroup of a finite group $G$. We say that $H$ satisfies partial $\Pi$-property in $G$ if there exists a chief series $\mathit{\Gamma}_G:1=G_0<G_1<\cdots<G_n=G$ of $G$ such that for every $G$-chief factor $G_i/G_{i-1}$ ($1\leq i\leq n$) of $\mathit{\Gamma}_G$, $|G/G_{i-1}:N_{G/G_{i-1}}(HG_{i-1}/G_{i-1}\cap G_i/G_{i-1})|$ is a $\pi(HG_{i-1}/G_{i-1}\cap G_i/G_{i-1})$-number. Our main results are listed here:\par
\medskip
\noindent\textbf{Theorem A. Let $\mathfrak{F}$ be a solubly saturated formation containing $\mathfrak{U}$ and $E$ a normal subgroup of $G$ with $G/E\in \mathfrak{F}$. Let $X\unlhd G$ such that $F_p^*(E)\leq X\leq E$. Suppose that for any Sylow $p$-subgroup $P$ of $X$, every maximal subgroup of $P$ satisfies partial $\Pi$-property in $G$. Then one of the following holds:}\par
\textbf{(1) $G\in \mathfrak{G}_{p'}\mathfrak{F}$.}\par
\textbf{(2) $X/O_{p'}(X)$ is a quasisimple group with Sylow $p$-subgroups of order $p$. In particular, if $X=F_p^*(E)$, then $X/O_{p'}(X)$ is a simple group.}\par
\medskip
\noindent\textbf{Theorem B. Let $\mathfrak{F}$ be a solubly saturated formation containing $\mathfrak{U}$ and $E$ a normal subgroup of $G$ with $G/E\in \mathfrak{F}$. Suppose that for any Sylow $p$-subgroup $P$ of $F_p^*(E)$, every cyclic subgroup of $P$ of prime order or order 4 (when $P$ is not quaternion-free) satisfies partial $\Pi$-property in $G$. Then $G\in \mathfrak{G}_{p'}\mathfrak{F}$.}\par
\end{abstract}
\renewcommand{\thefootnote}{\empty}

\textcolor{red}{\textbf{This manuscript has been published in J. Group Theory ([J. Group Theory 16 (2013), 745-766]). Now I am very regretful to say that there exists several mistakes in this paper. The following is a corrected and improved version.}}\par
\section{Introduction}
\noindent All groups considered in this paper are finite, $G$ always denotes a group and $p$ denotes a prime. Let $\pi$ denote a set of some primes and $\pi(G)$ denote the set of all prime divisors of $|G|$. $G_p$ denotes a Sylow $p$-subgroup of $G$ and $|G|_p$ denotes the order of $G_p$. An integer $n$ is called a $\pi$-number if all prime divisors of $n$ belong to $\pi$.\par
Recall that a class of groups $\mathfrak{F}$ is called a formation if $\mathfrak{F}$ is closed under taking homomorphic image and subdirect product. A formation $\mathfrak{F}$ is said to be saturated (resp. solubly saturated) if $G\in \mathfrak{F}$ whenever $G/\Phi(G)\in \mathfrak{F}$ (resp. $G/\Phi(N)\in \mathfrak{F}$ for a solvable normal subgroup $N$ of $G$). A $G$-chief factor $L/K$ is said to be $\mathfrak{F}$-central (resp. $\mathfrak{F}$-eccentric) in $G$ if $(L/K)\rtimes (G/C_G(L/K))\in \mathfrak{F}$ (resp. $(L/K)\rtimes (G/C_G(L/K))\notin \mathfrak{F}$). Following \cite{Gu10}, a normal subgroup $N$ of $G$ is called $\pi\mathfrak{F}$-hypercentral in $G$ if every $G$-chief factor below $N$ of order divisible by at least one prime in $\pi$ is $\mathfrak{F}$-central in $G$. Let $Z_{\pi\mathfrak{F}}(G)$ denote the $\pi\mathfrak{F}$-hypercentre of $G$, that is, the product of all $\pi\mathfrak{F}$-hypercentral normal subgroups of $G$. Let $Z_\mathfrak{F}(G)$ denote the $\mathfrak{F}$-hypercentre of $G$, that is, $Z_\mathfrak{F}(G)=Z_{\mathbb{P}\mathfrak{F}}(G)$.\par
We use $\mathfrak{U}$ (resp. $\mathfrak{U}_p$) to denote the class of finite supersolvable (resp. $p$-supersolvable) groups and $\mathfrak{N}$ (resp. $\mathfrak{N}_p$) to denote the class of finite nilpotent (resp. $p$-nilpotent) groups. Also, the symbol $\mathfrak{G}_\pi$ denotes the class of all finite $\pi$-groups. All notations and terminology not mentioned are standard, as in \cite{Hup1,Doe,Guo1}.\par
In \cite{Li1}, Li introduced the concept of \textit{$\Pi$-property} as follows: a subgroup $H$ of $G$ is said to satisfy $\Pi$-property in $G$ if for every $G$-chief factor $L/K$, $|G/K:N_{G/K}(HK/K\cap L/K)|$ is a $\pi(HK/K\cap L/K)$-number. Now we introduce the following concept which generalizes a large number of known embedding property (see Section 7 below).\par
\medskip
\noindent\textbf{Definition 1.1.} A subgroup $H$ of $G$ is said to satisfy \textit{partial $\Pi$-property} in $G$ if there exists a chief series $$\mathit{\Gamma}_G:1=G_0<G_1<\cdots<G_n=G$$ of $G$ such that for every $G$-chief factor $G_i/G_{i-1}$ ($1\leq i\leq n$) of $\mathit{\Gamma}_G$, $|G/G_{i-1}:N_{G/G_{i-1}}(HG_{i-1}/G_{i-1}\cap G_i/G_{i-1})|$ is a $\pi(HG_{i-1}/G_{i-1}\cap G_i/G_{i-1})$-number.\par
\medskip
Obviously, a subgroup $H$ of $G$ which satisfies $\Pi$-property in $G$ also satisfies partial $\Pi$-property in $G$. However, the converse does not hold as the following example illustrates.\par
\medskip
\noindent\textbf{Example 1.2.} Let $L_1=\langle a,b \,|\,a^5=b^5=1,ab=ba \rangle$ and $L_2=\langle a',b' \rangle$ be a copy of $L_1$. Let $\alpha$ be an automorphism of $L_1$ of order 3 satisfying that $a^\alpha=b$, $b^\alpha=a^{-1}b^{-1}$. Put $G=(L_1\times L_2)\rtimes \langle \alpha \rangle$. For any subgroup $H$ of $G$ of order 25, there exists a minimal normal subgroup $N$ of $G$ such that $H\cap N=1$ (for details, see \cite[Example]{Guo5}). Note that $\mathit{\Gamma}_G:1<N<HN<G$ is a chief series of $G$. Then $H$ satisfies partial $\Pi$-property in $G$. Now let $H'=\langle a \rangle\times \langle a' \rangle$. Since $|G:N_G(H'\cap L_1)|=|G:N_G(\langle a \rangle)|=3$, we have that $H'$ does not satisfy $\Pi$-property in $G$.\par
\medskip
Recall that the generalized Fitting subgroup $F^*(G)$ of $G$ is the quasinilpotent radical of $G$ (for details, see \cite[Chapter X]{Hup}). Note that $G$ is said to be $p$-quasinilpotent if $G$ induces inner automorphisms on each of its chief factors of order divisible by $p$. Following \cite{Laf}, the $p$-generalized Fitting subgroup $F_p^*(G)$ of $G$ is the $p$-quasinilpotent radical of $G$.\par
In this paper, we arrive at the following main results.\par
\medskip
\noindent\textbf{Theorem A.} Let $\mathfrak{F}$ be a solubly saturated formation containing $\mathfrak{U}$ and $E$ a normal subgroup of $G$ with $G/E\in \mathfrak{F}$. Let $X\unlhd G$ such that $F_p^*(E)\leq X\leq E$. Suppose that for any Sylow $p$-subgroup $P$ of $X$, every maximal subgroup of $P$ satisfies partial $\Pi$-property in $G$. Then one of the following holds:\par
(1) $G\in \mathfrak{G}_{p'}\mathfrak{F}$.\par
(2) $X/O_{p'}(X)$ is a quasisimple group with Sylow $p$-subgroups of order $p$. In particular, if $X=F_p^*(E)$, then $X/O_{p'}(X)$ is a simple group.\par
\medskip
\noindent\textbf{Theorem B.} Let $\mathfrak{F}$ be a solubly saturated formation containing $\mathfrak{U}$ and $E$ a normal subgroup of $G$ with $G/E\in \mathfrak{F}$. Suppose that for any Sylow $p$-subgroup $P$ of $F_p^*(E)$, every cyclic subgroup of $P$ of prime order or order 4 (when $P$ is not quaternion-free) satisfies partial $\Pi$-property in $G$. Then $G\in \mathfrak{G}_{p'}\mathfrak{F}$.\par
\medskip
\noindent\textbf{Theorem C.} Let $\mathfrak{F}$ be a solubly saturated formation containing $\mathfrak{U}$ and $E$ a normal subgroup of $G$ with $G/E\in \mathfrak{F}$. Let $X\unlhd G$ such that $F^*(E)\leq X\leq E$. Suppose that for any non-cyclic Sylow subgroup $P$ of $X$, either every maximal subgroup of $P$ satisfies partial $\Pi$-property in $G$, or every cyclic subgroup of $P$ of prime order or order 4 (when $P$ is not quaternion-free) satisfies partial $\Pi$-property in $G$. Then $G\in \mathfrak{F}$.\par
\medskip
The following propositions are the main stages of the proof of the above main results.\par
\medskip
\noindent\textbf{Proposition 1.3.} Suppose that $P$ is a normal $p$-subgroup of $G$. If every maximal subgroup of $P$ satisfies partial $\Pi$-property in $G$, then $P\leq Z_\mathfrak{U}(G)$.\par
\medskip
\noindent\textbf{Proposition 1.4.} Let $E$ be a normal subgroup of $G$ and $P$ a Sylow $p$-subgroup of $E$. If every maximal subgroup of $P$ satisfies partial $\Pi$-property in $G$, then either $E\leq Z_{p\mathfrak{U}}(G)$ or $|E|_p=p$.\par
\medskip
\noindent\textbf{Proposition 1.5.} Suppose that $P$ is a normal $p$-subgroup of $G$. If every cyclic subgroup of $P$ of prime order or order 4 (when $P$ is not quaternion-free) satisfies partial $\Pi$-property in $G$, then $P\leq Z_\mathfrak{U}(G)$.\par
\medskip
\noindent\textbf{Proposition 1.6.} Let $E$ be a normal subgroup of $G$ and $P$ a Sylow $p$-subgroup of $E$. If every cyclic subgroup of $P$ of prime order or order 4 (when $P$ is not quaternion-free) satisfies partial $\Pi$-property in $G$, then $E\leq Z_{p\mathfrak{U}}(G)$.\par
\medskip
\noindent\textbf{Proposition 1.7.} Let $E$ be a normal subgroup of $G$ and $P$ a Sylow $p$-subgroup of $E$ with $(|E|,p-1)=1$. If either every maximal subgroup of $P$ satisfies partial $\Pi$-property in $G$, or every cyclic subgroup of $P$ of prime order or order 4 (when $P$ is not quaternion-free) satisfies partial $\Pi$-property in $G$, then $E\in \mathfrak{N}_p$.\par
\medskip
Finally, we list the following corollaries which can be deduced from our theorems.\par
\medskip
\noindent\textbf{Corollary 1.8.} Let $\mathfrak{F}$ be a formation containing $\mathfrak{N}_p$ which satisfies $\mathfrak{G}_{p'}\mathfrak{F}=\mathfrak{F}$ and $E$ a normal subgroup of $G$ with $G/E\in \mathfrak{F}$. Suppose that for any Sylow $p$-subgroup $P$ of $E$, $N_G(P)\in \mathfrak{N}_p$ and either every maximal subgroup of $P$ satisfies partial $\Pi$-property in $G$, or every cyclic subgroup of $P$ of prime order or order 4 (when $P$ is not quaternion-free) satisfies partial $\Pi$-property in $G$. Then $G\in \mathfrak{F}$.\par
\medskip
\noindent\textbf{Corollary 1.9.} Let $\mathfrak{F}$ be a solubly saturated formation containing $\mathfrak{N}$ and $E$ a normal subgroup of $G$ with $G/E\in \mathfrak{F}$. Suppose that every subgroup of $F^*(E)$ of prime order is contained in $Z_\infty(G)$ and every cyclic subgroup of $F^*(E)$ of order 4 (when the Sylow 2-subgroups of $F^*(E)$ are not quaternion-free) satisfies partial $\Pi$-property in $G$. Then $G\in \mathfrak{F}$.\par
\section{Preliminaries}
\noindent\textbf{Lemma 2.1.} Let $H\leq G$ and $N\unlhd G$. Then:\par
(1) If $H$ satisfies partial $\Pi$-property in $G$, then $H^g$ satisfies partial $\Pi$-property in $G$ for every element $g\in G$.\par
(2) If $H$ is a $p$-subgroup of $G$, $H\leq N$ and $H$ satisfies partial $\Pi$-property in $G$, then $H$ satisfies partial $\Pi$-property in $N$.\par
(3) If either $N\leq H$ or $(|H|,|N|)=1$ and $H$ satisfies partial $\Pi$-property in $G$, then $HN/N$ satisfies partial $\Pi$-property in $G/N$.\par
(4) Let $P$ be a Sylow $p$-subgroup of $H$. If every maximal subgroup of $P$ satisfies partial $\Pi$-property in $G$ and $N\leq H$, then every maximal subgroup of $PN/N$ satisfies partial $\Pi$-property in $G/N$.\par
(5) If $HN/N$ satisfies partial $\Pi$-property in $G/N$ and there exists a chief series $\mathit{\Gamma}_N: 1=N_0<N_1<\cdots<N_n=N$ of $G$ below $N$ such that for every $G$-chief factor $N_i/N_{i-1}$ ($1\leq i\leq n$) of $\mathit{\Gamma}_N$, $|G:N_G(HN_{i-1}\cap N_i)|$ is a $\pi((HN_{i-1}\cap N_i)/N_{i-1})$-number, then $H$ satisfies partial $\Pi$-property in $G$.\par
\noindent\textit{Proof.} Statements (1) and (5) are obvious.\par
(2) Suppose that $H$ is a $p$-subgroup of $G$, $H\leq N$ and $H$ satisfies partial $\Pi$-property in $G$. Then there exists a chief series $$\mathit{\Gamma}_G:1=G_0<G_1<\cdots<G_n=G$$ of $G$ such that for every $G$-chief factor $G_i/G_{i-1}$ ($1\leq i\leq n$) of $\mathit{\Gamma}_G$, $|G:N_G(HG_{i-1}\cap G_i)|$ is a $p$-number. Now consider the normal series $$\mathit{\Gamma}_N:1=G_0\cap N\leq G_1\cap N\leq \cdots\leq G_n\cap N=N$$ of $N$. Avoiding repetitions, for every normal section $(G_i\cap N)/(G_{i-1}\cap N)$ ($1\leq i\leq n$) of $\mathit{\Gamma}_N$, we have that $H(G_{i-1}\cap N)\cap (G_i\cap N)=(H\cap G_i)(G_{i-1}\cap N)=(H\cap G_i)G_{i-1}\cap N$. It follows that $N_G((H\cap G_i)G_{i-1})=N_G((H\cap G_i)G_{i-1}\cap N)$, and so $|N:N_N(H(G_{i-1}\cap N)\cap (G_i\cap N))|$ is a $p$-number. Let $L/K$ be an $N$-chief factor such that $G_{i-1}\cap N\leq K\leq L\leq G_{i}\cap N$. Note that $$N_N(H(G_{i-1}\cap N)\cap (G_i\cap N))\leq N_N((H(G_{i-1}\cap N)\cap L)K)=N_N(HK\cap L).$$ Therefore, we get that $|N:N_N(HK\cap L)|$ is a $p$-number. This shows that $H$ satisfies partial $\Pi$-property in $N$.\par
(3) Suppose that either $N\leq H$ or $(|H|,|N|)=1$ and $H$ satisfies partial $\Pi$-property in $G$. Then for every normal subgroup $X$ of $G$, we have that $HN\cap XN=(H\cap X)N$. As $H$ satisfies partial $\Pi$-property in $G$, there exists a chief series $$\mathit{\Gamma}_G:1=G_0<G_1<\cdots<G_n=G$$ of $G$ such that for every $G$-chief factor $G_i/G_{i-1}$ ($1\leq i\leq n$) of $\mathit{\Gamma}_G$, $|G:N_G(HG_{i-1}\cap G_i)|$ is a $\pi((HG_{i-1}\cap G_i)/G_{i-1})$-number. Now consider the normal series $$\mathit{\Gamma}_{G/N}:1=G_0N/N\leq G_1N/N\leq \cdots\leq G_nN/N=G/N$$ of $G/N$. Avoiding repetitions, for every normal section $(G_iN/N)/(G_{i-1}N/N)$ ($1\leq i\leq n$) of $\mathit{\Gamma}_{G/N}$, we have that $HG_{i-1}N\cap G_iN=(H\cap G_i)G_{i-1}N$. Obviously, $N_G((H\cap G_i)G_{i-1})\leq N_G((H\cap G_i)G_{i-1}N)$. Hence $|G:N_G(HG_{i-1}N\cap G_iN)|$ is a $\pi((H\cap G_i)G_{i-1}/G_{i-1})$-number. Note that $G_i\cap G_{i-1}N=G_{i-1}$ for $G_iN\neq G_{i-1}N$. Then it is easy to see that $$\pi((HG_{i-1}N\cap G_iN)/G_{i-1}N)=\pi(H\cap G_i/H\cap G_i\cap G_{i-1}N)=\pi((H\cap G_i)G_{i-1}/G_{i-1}).$$ This shows that $HN/N$ satisfies partial $\Pi$-property in $G/N$.\par
(4) Let $T/N$ be a maximal subgroup of $PN/N$. Then $T/N=P_1N/N$, where $P_1$ is a maximal subgroup of $P$ such that $P_1\cap N=P\cap N$. By hypothesis, $P_1$ satisfies partial $\Pi$-property in $G$. Note that $P_1N\cap XN=(P_1\cap X)N$ for every normal subgroup $X$ of $G$. Similarly as the proof of (3), we can obtain that $T/N=P_1N/N$ satisfies partial $\Pi$-property in $G/N$, and thus (4) holds.\par
\medskip
If $P$ is either an odd order $p$-group or a quaternion-free 2-group, then let ${\Omega}(P)$ denote the subgroup ${\Omega}_1(P)$, otherwise ${\Omega}(P)$ denotes ${\Omega}_2(P)$.\par
\medskip
\noindent\textbf{Lemma 2.2.} \cite[Lemma 2.8]{Che3} Let $\mathfrak{F}$ be a solubly saturated formation, $P$ a normal $p$-subgroup of $G$ and $C$ a Thompson critical subgroup of $P$. If either $P/\Phi(P)\leq Z_\mathfrak{F}(G/\Phi(P))$ or ${\Omega}(C)\leq Z_\mathfrak{F}(G)$, then $P\leq Z_\mathfrak{F}(G)$.\par
\medskip
\noindent\textbf{Lemma 2.3.} \cite[Lemma 2.1.6]{Bal4} Let $G$ be a $p$-supersolvable group. Then $G'$ is $p$-nilpotent. In particular, if $O_{p'}(G)=1$, then $G$ has a unique Sylow $p$-subgroup.\par
\medskip
\noindent\textbf{Lemma 2.4.} \cite[VI, Theorem 14.3]{Hup1} Suppose that $G$ has an abelian Sylow $p$-subgroup $P$. Then $G'\cap Z(G)\cap P=1$.\par
\medskip
\noindent\textbf{Lemma 2.5.} Let $\mathfrak{F}$ be any formation and $E\unlhd G$.\par
(1) \cite[Theorem B]{Ski} If $F^*(E)\leq Z_\mathfrak{F}(G)$, then $E\leq Z_\mathfrak{F}(G)$.\par
(2) \cite[Lemma 2.13]{Su} If $F_p^*(E)\leq Z_{p\mathfrak{F}}(G)$, then $E\leq Z_{p\mathfrak{F}}(G)$.\par
\medskip
\noindent\textbf{Lemma 2.6.} \cite[Lemma 2.11]{Su} Let $\mathfrak{F}$ be a solubly saturated formation containing $\mathfrak{U}$. Suppose that $E\unlhd G$ with $G/E\in \mathfrak{F}$.\par
(1) If $E\leq Z_\mathfrak{U}(G)$, then $G\in \mathfrak{F}$.\par
(2) If $E\leq Z_{p\mathfrak{U}}(G)$, then $G\in \mathfrak{G}_{p'}\mathfrak{F}$.\par
\medskip
\noindent\textbf{Lemma 2.7.} \cite[Lemma 3.1]{Bal6} Let $G$ be a group whose Sylow $p$-subgroups are cyclic groups of order $p$. Then: either (1) $G$ is a $p$-solvable group, or (2) $G/O_{p'}(G)$ is a quasisimple group such that $Soc(G/O_{p'}(G))=O^{p'}(G/O_{p'}(G))$ is a simple group whose Sylow $p$-subgroups are cyclic groups of order $p$.\par
\medskip
\noindent\textbf{Lemma 2.8.} \cite[Lemma 3.1]{War} Let $G$ be a non-abelian quaternion-free $2$-group. Then $G$ has a characteristic subgroup of index $2$.\par
\medskip
\noindent\textbf{Lemma 2.9.} \cite[Lemma 2.10]{Che3} Let $C$ be a Thompson critical subgroup of a nontrivial $p$-group $P$.\par
(1) If $p$ is odd, then the exponent of ${\Omega}_1(C)$ is $p$.\par
(2) If $P$ is an abelian $2$-group, then the exponent of ${\Omega}_1(C)$ is $2$.\par
(3) If $p=2$, then the exponent of ${\Omega}_2(C)$ is at most $4$.\par
\medskip
Following \cite{Bal5}, let $\Psi_p(G)=\langle x\mid x\in G, o(x)=p\rangle$ if $p$ is odd, and $\Psi_2(G)=\langle x\mid x\in G, o(x)=2$ or $4\rangle$.\par
\noindent\textbf{Lemma 2.10.} \cite[Theorem 6]{Bal5} Let $K$ be a normal subgroup of $G$ with $G/K$ contained in a saturated formation $\mathfrak{F}$. If $\Psi_p(K)\leq Z_\mathfrak{F}(G)$, then $G/O_{p'}(K)\in \mathfrak{F}$.\par
\medskip
\noindent\textbf{Lemma 2.11.} \cite[Corollary 2]{Asa} Let $P$ be a Sylow 2-subgroup of $G$. If $P$ is quaternion-free and ${\Omega}_1(P)\leq Z(G)$, then $G$ is 2-nilpotent.\par
\medskip
\noindent\textbf{Lemma 2.12.} Let $p$ be a prime divisor of $|G|$ with $(|G|,p-1)=1$. If $G$ has cyclic Sylow $p$-subgroups, then $G$ is $p$-nilpotent.\par
\noindent\textit{Proof.} Prove similarly as \cite[(10.1.9)]{Rob}.\par
\section{Proof of Theorem A}
\noindent\textbf{Proof of Proposition 1.3.} Suppose that the result is false and let $(G,P)$ be a counterexample for which $|G|+|P|$ is minimal. We proceed via the following steps.\par
(1) $G$ has a unique minimal normal subgroup $N$ contained in $P$, $P/N\leq Z_\mathfrak{U}(G/N)$ and $|N|>p$.\par
Let $N$ be a minimal normal subgroup of $G$ contained in $P$. By Lemma 2.1(3), $(G/N,P/N)$ satisfies the hypothesis. Then $P/N\leq Z_\mathfrak{U}(G/N)$ by the choice of $(G,P)$. If $|N|=p$, then $N\leq Z_\mathfrak{U}(G)$, and so $P\leq Z_\mathfrak{U}(G)$, which is absurd. Hence $|N|>p$. Now suppose that $G$ has a minimal normal subgroup $R$ contained in $P$, which is different from $N$. Then $P/R\leq Z_\mathfrak{U}(G/R)$ as above. It follows that $NR/R\leq Z_\mathfrak{U}(G/R)$, and thereby $N\leq Z_\mathfrak{U}(G)$ for $G$-isomorphism $NR/R\cong N$. Therefore, we have that $P\leq Z_\mathfrak{U}(G)$, a contradiction. This shows that (1) holds.\par
(2) $\Phi(P)\neq 1$.\par
If $\Phi(P)=1$, then $P$ is elementary abelian. This induces that $N$ has a complement $S$ in $P$. Let $L$ be a maximal subgroup of $N$ such that $L$ is normal in some Sylow $p$-subgroup $G_p$ of $G$. Then $L\neq 1$ and $H=LS$ is a maximal subgroup of $P$. By hypothesis, $H$ satisfies partial $\Pi$-property in $G$. Then $G$ has a chief series $\mathit{\Gamma}_G:1=G_0<G_1<\cdots<G_n=G$ such that for every $G$-chief factor $G_i/G_{i-1}$ ($1\leq i\leq n$) of $\mathit{\Gamma}_G$, $|G:N_G(HG_{i-1}\cap G_i)|$ is a $p$-number. Note that there exists an integer $k$ ($1\leq k\leq n$) such that $G_k=G_{k-1}\times N$. It follows that $|G:N_G(HG_{k-1}\cap G_k)|$ is a $p$-number, and so $|G:N_G(HG_{k-1}\cap N)|$ is a $p$-number. Since $L\leq HG_{k-1}\cap N\leq N$, we have that either $HG_{k-1}\cap N=N$ or $HG_{k-1}\cap N=L$. In the former case, if $G_{k-1}\cap P\neq 1$, then $N\leq G_{k-1}$ by (1), which is impossible. Hence $G_{k-1}\cap P=1$, and thus $N\leq H(G_{k-1}\cap P)=H$, a contradiction. In the latter case, since $L\unlhd G_p$, we get that $L\unlhd G$, also a contradiction. Thus (2) follows.\par
(3) Final contradiction.\par
Since $\Phi(P)\neq 1$ and $N$ is the unique minimal normal subgroup of $G$ contained in $P$, we have that $N\leq \Phi(P)$. This deduces that $P/\Phi(P)\leq Z_\mathfrak{U}(G/\Phi(P))$. Then by Lemma 2.2, $P\leq Z_\mathfrak{U}(G)$. The final contradiction completes the proof.\par
\medskip
\noindent\textbf{Proof of Proposition 1.4.} Suppose that the result is false and let $(G,E)$ be a counterexample for which $|G|+|E|$ is minimal. We proceed via the following steps.\par
(1) $O_{p'}(E)=1$.\par
If $O_{p'}(E)\neq 1$, then the hypothesis holds for $(G/O_{p'}(E),E/O_{p'}(E))$ by Lemma 2.1(3). The choice of $(G,E)$ implies that either $E/O_{p'}(E)\leq Z_{p\mathfrak{U}}(G/O_{p'}(E))=Z_{p\mathfrak{U}}(G)/O_{p'}(E)$ or $|E/O_{p'}(E)|_p=p$. It follows that either $E\leq Z_{p\mathfrak{U}}(G)$ or $|E|_p=p$, a contradiction.\par
(2) $E=G$.\par
Suppose that $E<G$. By Lemma 2.1(2), $(E,E)$ satisfies the hypothesis. By the choice of $(G,E)$, either $E\in \mathfrak{U}_p$ or $|E|_p=p$. We may, therefore, assume that $E\in \mathfrak{U}_p$. Then by (1) and Lemma 2.3, we get that $P\unlhd E$, and so $P\unlhd G$. By Proposition 1.3, we have that $P\leq Z_\mathfrak{U}(G)$. This induces that $E\leq Z_{p\mathfrak{U}}(G)$, which is absurd.\par
(3) $G$ has a unique minimal normal subgroup $N$, $p\mid|N|$, and either $G/N\in \mathfrak{U}_p$ or $|G/N|_p=p$.\par
Let $N$ be a minimal normal subgroup of $G$. By Lemma 2.1(4), $(G/N,G/N)$ satisfies the hypothesis. By the choice of $(G,E)$, we have that either $G/N\in \mathfrak{U}_p$ or $|G/N|_p=p$. Since $O_{p'}(G)=1$, we get that $p\mid |N|$. Let $R$ be a minimal normal subgroup of $G$, which is different from $N$. Then $p\mid|R|$ and either $G/R\in \mathfrak{U}_p$ or $|G/R|_p=p$ as above. First suppose that $G/N\in \mathfrak{U}_p$ and $G/R\in \mathfrak{U}_p$. Then $G\in \mathfrak{U}_p$, a contradiction.\par
Next consider that $G/N\in \mathfrak{U}_p$ and $|G/R|_p=p$. Note that $RN/N$ is a minimal normal subgroup of $G/N$ and $p\mid|R|$. This induces that $|R|=|RN/N|=p$, and so $|P|=|G|_p=p^2$. Since $|N|_p=|NR/R|_p\leq |G/R|_p=p$ and $p\mid |N|$, we have that $|N|_p=p$. This shows that $P\cap N\in Syl_p(N)$ is a nontrivial maximal subgroup of $P$. Hence $P\cap N$ satisfies partial $\Pi$-property in $G$. Then $G$ has a chief series $\mathit{\Gamma}_G:1=G_0<G_1<\cdots<G_n=G$ such that for every $G$-chief factor $G_i/G_{i-1}$ ($1\leq i\leq n$) of $\mathit{\Gamma}_G$, $|G:N_G((P\cap N)G_{i-1}\cap G_i)|$ is a $p$-number. Note that there exists an integer $k$ ($1\leq k\leq n$) such that $G_k=G_{k-1}\times N$. It follows that $|G:N_G((P\cap N)G_{k-1}\cap G_k)|$ is a $p$-number. Since $P\cap N\unlhd P$, we get that $(P\cap N)G_{k-1}\cap G_k\unlhd G$. This deduces that $P\cap N\unlhd G$, and so $N\leq P$. Consequently, $|N|=p$. It follows that $G\in \mathfrak{U}_p$, which is absurd. If $G/R\in \mathfrak{U}_p$ and $|G/N|_p=p$, we can handle it in a similar way.\par
Finally, assume that $|G/N|_p=p$ and $|G/R|_p=p$. Then since $p\mid |N|$ and $p\mid |R|$, we get that $|N|_p=|R|_p=p$ and $|G|_p=p^2$. This induces that $P\cap N$ and $P\cap R$ are nontrivial maximal subgroups of $P$, and so $P\cap N$ and $P\cap R$ satisfy partial $\Pi$-property in $G$. With a similar discussion as above, we have that $P\cap N\unlhd G$ and $P\cap R\unlhd G$. This implies that $N\leq P$ and $R\leq P$. Therefore, $P=N\times R\unlhd G$, and thereby $G\in \mathfrak{U}_p$. The final contradiction shows that (3) holds.\par
(4) $N\leq O_p(G)$.\par
If not, then $O_p(G)=1$ by (3). Let $H$ be a maximal subgroup of $P$. Then $H$ satisfies partial $\Pi$-property in $G$. Thus $G$ has a chief series $\mathit{\Gamma}_G:1=G_0<G_1=N<\cdots<G_n=G$ such that for every $G$-chief factor $G_i/G_{i-1}$ ($1\leq i\leq n$) of $\mathit{\Gamma}_G$, $|G:N_G(HG_{i-1}\cap G_i)|$ is a $p$-number. It follows that $|G:N_G(H\cap N)|$ is a $p$-number. Since $H\cap N\unlhd P$, we have that $H\cap N\unlhd G$, and so $H\cap N\leq O_p(G)=1$. Hence $|N|_p=p$. If $|G|_p>p$, then $P$ has a maximal subgroup $L$ containing $P\cap N$. By hypothesis, $L$ satisfies partial $\Pi$-property in $G$. Similarly as above, we can conclude that $L\cap N\leq O_p(G)=1$, and thereby $P\cap N=L\cap N=1$, a contradiction. Therefore, $|G|_p=p$, which is absurd. Thus (4) follows.\par
(5) $N\nleq \Phi(P)$.\par
If $N\leq \Phi(P)$, then $N\leq \Phi(G)$. If $G/N\in \mathfrak{U}_p$, then $G\in \mathfrak{U}_p$, which is impossible. Hence by (3), we may assume that $|G/N|_p=p$. Put $A/N=O_{p'}(G/N)$. Since $A\cap P\leq N\leq \Phi(P)$, $A$ is $p$-nilpotent by \cite[IV, Theorem 4.7]{Hup1}. Let $A_{p'}$ be the normal $p$-complement of $A$. Then $A_{p'}\leq O_{p'}(G)=1$, and thus $A$ is a $p$-group. It follows that $A=N$, and so $O_{p'}(G/N)=1$. Let $X/N$ be a $G$-chief factor. As $O_{p'}(G/N)=1$, we have that $p\mid |X/N|$. This deduces that $|X/N|_p=|G/N|_p=p$ and $P\leq X$. Clearly, the hypothesis holds for $(G,X)$. Suppose that $X<G$. Then by the choice of $(G,E)$, we obtain that either $X\leq Z_{p\mathfrak{U}}(G)$ or $|X|_p=p$. In the former case, $G\in \mathfrak{U}_p$, a contradiction. In the latter case, $|G|_p=|X|_p=p$, also a contradiction. Thus $X=G$. Then $G/N$ is a $G$-chief factor. Considering the above, we may assume that $G/N$ is a non-abelian simple group.\par
It is clear that $N$ is a maximal subgroup of $P$. Thus $N=\Phi(P)$, and so $P$ is a cyclic group. It follows that $|N|=p$ and $|G|_p=|P|=p^2$. As $G/N$ is a non-abelian simple group, we have that $G'=G$. Note that $G/C_G(N)\apprle Aut(N)$ is abelian. This implies that $C_G(N)=G'=G$, and so $N\leq Z(G)$. Hence $N\leq G'\cap Z(G)\cap P$, which contradicts Lemma 2.4. This ends the proof of (5).\par
(6) Final contradiction.\par
Since $N\nleq \Phi(P)$, $P$ has a maximal subgroup $H$ such that $N\nleq H$. By hypothesis, $H$ satisfies partial $\Pi$-property in $G$. Thus $G$ has a chief series $\mathit{\Gamma}_G:1=G_0<G_1=N<\cdots<G_n=G$ such that for every $G$-chief factor $G_i/G_{i-1}$ ($1\leq i\leq n$) of $\mathit{\Gamma}_G$, $|G:N_G(HG_{i-1}\cap G_i)|$ is a $p$-number. It follows that $|G:N_G(H\cap N)|$ is a $p$-number. As $H\cap N\unlhd P$, we get that $H\cap N\unlhd G$. Therefore, $H\cap N=1$, and so $|N|=p$.\par
If $G/N\in \mathfrak{U}_p$, then $G\in \mathfrak{U}_p$, a contradiction. Thus by (3), $|G/N|_p=p$ holds. Then there exists an integer $k$ ($2\leq k\leq n$) such that $p\mid |G_k/G_{k-1}|$. Without loss of generality, we may assume that $H\leq G_k$ and $G_k/G_{k-1}$ is a non-abelian simple group. By hypothesis, $|G:N_G(HG_{k-1})|$ is a $p$-number. Since $H\unlhd P$, we have that $HG_{k-1}\unlhd G$. It follows that either $HG_{k-1}=G_k$ or $HG_{k-1}=G_{k-1}$. In the former case, $|G_k/G_{k-1}|$ is a $p$-number, a contradiction. In the latter case, $H\leq G_{k-1}$, and so $P=HN\leq G_{k-1}$, also a contradiction. The proof is thus completed.\par
\medskip
\noindent\textbf{Proof of Theorem A.} By Proposition 1.4, we have that either $X\leq Z_{p\mathfrak{U}}(G)$ or $|X|_p=p$. If $X\leq Z_{p\mathfrak{U}}(G)$, then $F_p^*(E)\leq Z_{p\mathfrak{U}}(G)$. Hence by Lemma 2.5(2), $E\leq Z_{p\mathfrak{U}}(G)$, and so $G\in \mathfrak{G}_{p'}\mathfrak{F}$ by Lemma 2.6(2). Now consider that $|X|_p=p$. We may suppose that $X$ is not $p$-solvable. Then by Lemma 2.7, $X/O_{p'}(X)$ is a quasisimple group.\par
In additional, assume that $X=F_p^*(E)$. By \cite[Lemma 2.10(2)]{Bal6}, $X/O_{p'}(X)=F_p^*(X)/O_{p'}(X)=F_p^*(X/O_{p'}(X))=F^*(X/O_{p'}(X))$ is quasinilpotent. Since $X/O_{p'}(X)$ is not $p$-solvable, $Soc(X/O_{p'}(X))$ is a non-abelian simple group by Lemma 2.7, and so $F(X/O_{p'}(X))=1$. It follows from \cite[X, Theorem 13.13]{Hup} that $X/O_{p'}(X)=Soc(X/O_{p'}(X))$ is simple. Thus the theorem is proved.\par

\section{Proof of Theorem B}
\noindent\textbf{Proof of Proposition 1.5.} Suppose that the result is false and let $(G,P)$ be a counterexample for which $|G|+|P|$ is minimal. We proceed via the following steps.\par
(1) $G$ has a unique normal subgroup $N$ such that $P/N$ is a $G$-chief factor, $N\leq Z_\mathfrak{U}(G)$ and $|P/N|>p$.\par
Let $P/N$ be a $G$-chief factor. It is easy to see that $(G,N)$ satisfies the hypothesis. By the choice of $(G,P)$, we have that $N\leq Z_\mathfrak{U}(G)$. If $|P/N|=p$, then $P/N\leq Z_\mathfrak{U}(G/N)$, and so $P\leq Z_\mathfrak{U}(G)$, which is contrary to our assumption. Hence $|P/N|>p$. Now assume that $P/R$ is a $G$-chief factor, which is different from $P/N$. Then $R\leq Z_\mathfrak{U}(G)$ as above. By $G$-isomorphism $P/N=NR/N\cong R/N\cap R$, we have that $P/N\leq Z_\mathfrak{U}(G/N)$, a contradiction. Thus (1) follows.\par
(2) The exponent of $P$ is $p$ or 4 (when $P$ is not quaternion-free).\par
Let $C$ be a Thompson critical subgroup of $P$. If ${\Omega}(C)<P$, then ${\Omega}(C)\leq N\leq Z_\mathfrak{U}(G)$ by (1), and so $P\leq Z_\mathfrak{U}(G)$ by Lemma 2.2, which is impossible. Hence $P=C={\Omega}(C)$.
If $P$ is a non-abelian quaternion-free 2-group, then $P$ has a characteristic subgroup $T$ of index 2 by Lemma 2.8. By (1), $T\leq N$, and so $|P/N|=2$, which is absurd. Thus $P$ is a non-abelian 2-group if and only if $P$ is not quaternion-free. Then by Lemma 2.9, the exponent of $P$ is $p$ or 4 (when $P$ is not quaternion-free).\par
(3) Final contradiction.\par
Note that $P/N\cap Z(G_p/N)>1$. Let $L/N$ be a subgroup of $P/N\cap Z(G_p/N)$ of order $p$. Then we may choose an element $l\in L\backslash N$. Put $H=\langle l\rangle$. Then $L=HN$ and $H$ is a subgroup of order $p$ or 4 (when $P$ is not quaternion-free) by (2). By hypothesis, $H$ satisfies partial $\Pi$-property in $G$. Then $G$ has a chief series $\mathit{\Gamma}_G:1=G_0<G_1<\cdots<G_n=G$ such that for every $G$-chief factor $G_i/G_{i-1}$ ($1\leq i\leq n$) of $\mathit{\Gamma}_G$, $|G:N_G(HG_{i-1}\cap G_i)|$ is a $p$-number. Clearly, there exists an integer $k$ ($1\leq k\leq n$) such that $P\nleq G_{k-1}N$ and $P\leq G_{k}N$. Since $N$ is the unique normal subgroup of $G$ such that $P/N$ is a $G$-chief factor, we have that $P\cap G_{k-1}\leq N$ and $P\cap G_{k}=P$. Thus $|G:N_G(HG_{k-1})|$ is a $p$-number. As $H\unlhd G_p$, we obtain that $HG_{k-1}\unlhd G$. This implies that $L=HN=(HG_{k-1}\cap P)N\unlhd G$. Hence $|P/N|=|L/N|=p$, a contradiction. The proof is thus completed.\par
\medskip
\noindent\textbf{Proof of Proposition 1.6.} Suppose that the result is false and let $(G,E)$ be a counterexample for which $|G|+|E|$ is minimal. We proceed via the following steps.\par
(1) $O_{p'}(E)=1$.\par
Assume that $O_{p'}(E)\neq 1$. Then the hypothesis holds for $(G/O_{p'}(E),E/O_{p'}(E))$ by Lemma 2.1(3). By the choice of $(G,E)$, we have that $E/O_{p'}(E)\leq Z_{p\mathfrak{U}}(G/O_{p'}(E))=Z_{p\mathfrak{U}}(G)/O_{p'}(E)$. This implies that $E\leq Z_{p\mathfrak{U}}(G)$, a contradiction.\par
(2) $E=G$.\par
If $E<G$, then by Lemma 2.1(2), $(E,E)$ satisfies the hypothesis. Due to the choice of $(G,E)$, we get that $E\in \mathfrak{U}_p$. It follows from (1) and Lemma 2.3 that $P\unlhd G$. Then by Proposition 1.5, we have that $P\leq Z_\mathfrak{U}(G)$, and thereby $E\leq Z_{p\mathfrak{U}}(G)$, which is impossible.\par
(3) $O_p(G)\leq Z_\mathfrak{U}(G)$.\par
This follows directly from Proposition 1.5.\par
(4) $Z_{p\mathfrak{U}}(G)$ is the unique normal subgroup of $G$ such that $G/Z_{p\mathfrak{U}}(G)$ is a $G$-chief factor and $Z(G)=Z_{\mathfrak{U}}(G)=O_p(G)$ is the Sylow $p$-subgroup of $Z_{p\mathfrak{U}}(G)$.\par
Let $G/L$ be a $G$-chief factor. Then clearly, $(G,L)$ satisfies the hypothesis. By the choice of $(G,E)$, $L\leq Z_{p\mathfrak{U}}(G)$, and so $L=Z_{p\mathfrak{U}}(G)$ for $G\notin \mathfrak{U}_p$. This implies that $Z_{p\mathfrak{U}}(G)$ is the unique normal subgroup of $G$ such that $G/Z_{p\mathfrak{U}}(G)$ is a $G$-chief factor. Since $O_{p'}(G)=1$ by (1) and (2), $O_p(G)$ is the Sylow $p$-subgroup of $Z_{p\mathfrak{U}}(G)$ by (3) and Lemma 2.3. If $G^{\mathfrak{U}}<G$, then $G^{\mathfrak{U}}\leq Z_{p\mathfrak{U}}(G)$. It follows from Lemma 2.6(2) that $G\in \mathfrak{U}_p$, a contradiction. Thus $G^{\mathfrak{U}}=G$. By \cite[IV, Theorem 6.10]{Doe}, $Z_{\mathfrak{U}}(G)\leq Z(G)$, and thereby $Z_{\mathfrak{U}}(G)=Z(G)$. Since $O_{p'}(Z(G))\leq O_{p'}(G)=1$ and $O_p(G)\leq Z_{\mathfrak{U}}(G)$ by (1)-(3), we obtain that $Z(G)=O_p(G)$. Hence (4) holds.\par
(5) $\Psi_p(G)=G$.\par
If not, then since $\Psi_p(G)\unlhd G$, we have that $\Psi_p(G)\leq Z_{p\mathfrak{U}}(G)$ by (4), and so $\Psi_p(G)\leq Z_\mathfrak{U}(G)$. By Lemma 2.10, $G\in \mathfrak{U}_p$, which is absurd.\par
(6) Final contradiction.\par
First assume that either $p>2$ or $p=2$ and $P$ is not quaternion-free. Since $\Psi_p(G)=G$, there exists an element $x$ of $G$ of order $p$ or 4 not contained in $Z_{p\mathfrak{U}}(G)$. By Lemma 2.1(1), without loss of generality, we may let $x\in P$. Put $H=\langle x\rangle$. Then $H$ satisfies partial $\Pi$-property in $G$. Thus $G$ has a chief series $\mathit{\Gamma}_G:1=G_0<G_1<\cdots<G_{n-1}=Z_{p\mathfrak{U}}(G)<G_n=G$ such that for every $G$-chief factor $G_i/G_{i-1}$ ($1\leq i\leq n$) of $\mathit{\Gamma}_G$, $|G:N_G(HG_{i-1}\cap G_i)|$ is a $p$-number. It follows from (4) that $|G:N_G(HZ_{p\mathfrak{U}}(G))|$ is a $p$-number. This implies that $G=(HZ_{p\mathfrak{U}}(G))^G=(HZ_{p\mathfrak{U}}(G))^{P}\leq PZ_{p\mathfrak{U}}(G)$. Hence $|G/Z_{p\mathfrak{U}}(G)|=p$, and so $G=Z_{p\mathfrak{U}}(G)$, a contradiction.\par
Now consider that $p=2$ and $P$ is quaternion-free. Put ${{\Psi}_2}'(G)=\langle x\mid x\in G, o(x)=2\rangle$. Then since ${{\Psi}_2}'(G)\unlhd G$, either ${{\Psi}_2}'(G)\leq Z_{2\mathfrak{U}}(G)$ or ${{\Psi}_2}'(G)=G$ by (4). If ${{\Psi}_2}'(G)\leq Z_{2\mathfrak{U}}(G)$, then ${\Omega}_1(P)\leq Z_{\mathfrak{U}}(G)=Z(G)$ by (4), and so $G$ is 2-nilpotent by Lemma 2.11, which is impossible. Therefore, we have that ${\Psi'}_2(G)=G$. Then there exists an element $y$ of $G$ of order 2 not contained in $Z_{2\mathfrak{U}}(G)$. With a similar discussion as above, we can also obtain a contradiction. This completes the proof.\par
\medskip
\noindent\textbf{Proof of Theorem B.} By Proposition 1.6, we have that $F_p^*(E)\leq Z_{p\mathfrak{U}}(G)$. It follows from Lemma 2.5(2) that $E\leq Z_{p\mathfrak{U}}(G)$, and so the theorem holds by Lemma 2.6(2).\par
\section{Proof of Theorem C}
\noindent\textbf{Proof of Proposition 1.7.} By Lemma 2.12, Proposition 1.4, and Proposition 1.6, we have that $E\in \mathfrak{U}_p$. Since $(|E|,p-1)=1$, it is easy to see that every $E$-chief factor of order $p$ is central in $E$. Hence $E\in \mathfrak{N}_p$ holds.\par
\medskip
\noindent\textbf{Proof of Theorem C.} Suppose that the result is false and let $(G,E)$ be a counterexample for which $|G|+|E|$ is minimal. Let $p$ be the smallest prime divisor of $|X|$ and $P\in Syl_p(X)$. If $P$ is cyclic, then $X$ is $p$-nilpotent by Lemma 2.12. If $P$ is not cyclic, then by Proposition 1.7, $X$ is also $p$-nilpotent. Let $X_{p'}$ be the normal $p$-complement of $X$. Then $X_{p'}\unlhd G$. If $P$ is cyclic, then $X/X_{p'}\leq Z_\mathfrak{U}(G/X_{p'})$. Now consider that $P$ is not cyclic. Then by Lemma 2.1(3), $(G/X_{p'},X/X_{p'})$ satisfies the hypothesis of Proposition 1.3 or Proposition 1.5. Hence $X/X_{p'}\leq Z_\mathfrak{U}(G/X_{p'})$ also holds.\par
Let $q$ be the smallest prime divisor of $|X_{p'}|$ and $Q\in Syl_q(X)$. With a similar argument as above, we get that $X_{p'}$ is $q$-nilpotent and $X_{p'}/X_{\{p,q\}'}\leq Z_\mathfrak{U}(G/X_{\{p,q\}'})$, where $X_{\{p,q\}'}$ is the normal $q$-complement of $X_{p'}$. The rest may be deduced by analogy. Hence we obtain that $X\leq Z_\mathfrak{U}(G)$. It follows from Lemma 2.5(1) that $E\leq Z_\mathfrak{U}(G)$. Then by Lemma 2.6(1), $G\in \mathfrak{F}$, which completes the proof.\par
\section{Proof of the Corollaries}
\noindent\textbf{Proof of Corollary 1.8.} Suppose that the result is false and let $(G,E)$ be a counterexample for which $|G|+|E|$ is minimal. By Lemma 2.1(3), it is easy to see that the hypothesis holds for $(G/O_{p'}(E),E/O_{p'}(E))$. If $O_{p'}(E)\neq 1$, then by the choice of $(G,E)$, $G/O_{p'}(E)\in \mathfrak{F}$, and so $G\in \mathfrak{F}$, a contradiction. Hence $O_{p'}(E)=1$. By Proposition 1.4 and Proposition 1.6, we get that either $E\in \mathfrak{U}_p$ or $|E|_p=p$. Suppose that $|E|_p=p$. Then $P$ is a cyclic group of order $p$. Since $N_G(P)\in \mathfrak{N}_p$, we have that $N_E(P)=P\times H$, where $H$ is the normal $p$-complement of $N_E(P)$. It follows that $N_E(P)=C_E(P)$. Hence $E\in \mathfrak{N}_p$ by Burnside's Theorem. As $O_{p'}(E)=1$, $P=E\unlhd G$. This implies that $G=N_G(P)\in \mathfrak{N}_p\subseteq \mathfrak{F}$, a contradiction. We may, therefore, assume that $E\in \mathfrak{U}_p$. In this case, $P\unlhd E$ by Lemma 2.3, and so $P\unlhd G$. This induces that $G=N_G(P)\in \mathfrak{N}_p\subseteq \mathfrak{F}$, also a contradiction.\par
\medskip
\noindent\textbf{Proof of Corollary 1.9.} By \cite[Proposition 2.3(2)]{Li1}, every cyclic subgroup of $F^*(E)$ of prime order or order 4 (when the Sylow 2-subgroups of $F^*(E)$ are not quaternion-free) satisfies partial $\Pi$-property in $G$. It follows from Proposition 1.6 that $F^*(E)\leq Z_\mathfrak{U}(G)$, and so $F^*(E)=F(E)$ by \cite[X, Corollary 13.7(d)]{Hup}. Note that $O_2(E)\leq Z_\infty(G)$. Then by Lemma 2.2, we have that $F^*(E)=F(E)\leq Z_\infty(G)$. It follows from Lemma 2.5(1) that $E\leq Z_\infty(G)\leq Z_\mathfrak{F}(G)$. Therefore, the corollary holds by \cite[Lemma 2.13]{Guo}.\par
\section{Remarks and Applications}
In this section, we shall show that partial $\Pi$-property still holds on the subgroups which satisfy a certain known embedding property mentioned below. In brief, we only focus on most important and recent embedding properties.\par
Recall that a subgroup $H$ of $G$ is called to be a CAP-subgroup if $H$ either covers or avoids every $G$-chief factor. Let $\mathfrak{F}$ be a saturated formation. A subgroup $H$ of $G$ is said to be $\mathfrak{F}$-hypercentrally embedded \cite{Ezq} in $G$ if $H^G/H_G\leq Z_{\mathfrak{F}}(G/H_G)$. A subgroup $H$ of $G$ is called to be quasinormal (or permutable) in $G$ if $H$ $H$ permutes with every subgroup of $G$. A subgroup $H$ of $G$ is said to be S-quasinormal (or S-permutable) in $G$ if $H$ permutes with every Sylow subgroup of $G$. Let $X$ be a non-empty subset of $G$. A subgroup $H$ of $G$ is called to be $X$-permutable \cite{Guo6} with a subgroup $T$ of $G$ if there exists an element $x\in X$ such that $HT^x=T^xH$. A subgroup $H$ of $G$ is said to be S-semipermutable \cite{Che} in $G$ if $H$ permutes with every Sylow $p$-subgroup of $G$ such that $(p,|H|)=1$. A subgroup $H$ of $G$ is called to be SS-quasinormal \cite{Li} in $G$ if $H$ has a supplement $K$ in $G$ such that $H$ permutes with every Sylow subgroup of $K$.\par
\medskip
\noindent\textbf{Lemma 7.1.} Let $H$ be a subgroup of $G$. Then $H$ satisfies $\Pi$-property, and thus satisfies partial $\Pi$-property in $G$, if one of the following holds:\par
(1) $H$ is a CAP-subgroup of $G$.\par
(2) $H$ is $\mathfrak{U}$-hypercentrally embedded in $G$.\par
(3) $H$ is S-quasinormal in $G$.\par
(4) $H$ is $X$-permutable with all Sylow subgroups of $G$, where $X$ is a solvable normal subgroup of $G$.\par
(5) $H$ is a $p$-group and $H$ is S-semipermutable in $G$.\par
(6) $H$ is a $p$-group and $H$ is SS-quasinormal in $G$.\par
\noindent\textit{Proof.} Statements (1)-(4) were proved in \cite{Li1}, and the proof of \cite[Proposition 2.4]{Li1} still works for statement (5).\par
(6) By (5), we only need to prove that $H$ is S-semipermutable in $G$. By definition, $H$ has a supplement $K$ in $G$ such that $H$ permutes with every Sylow subgroup of $K$. Let $G_p$ be a Sylow $p$-subgroup of $G$ such that $(p,|H|)=1$. Then there exists an element $h\in H$ such that ${G_p}^h\leq K$. It follows that $H{G_p}^h={G_p}^hH$, and thereby $HG_p=G_pH$. Hence $H$ is S-semipermutable in $G$. This shows that $H$ satisfies $\Pi$-property in $G$.\par
\medskip
Recall that a subgroup $H$ of $G$ is called to be a partial CAP-subgroup (or semi CAP-subgroup) \cite{Fan} if there exists a chief series $\mathit{\Gamma}_G$ of $G$ such that $H$ either covers or avoids every $G$-chief factor of $\mathit{\Gamma}_G$. A subgroup $H$ of $G$ is said to be S-embedded \cite{Guo3} in $G$ if $G$ has a normal subgroup $K$ such that $HK$ is S-quasinormal in $G$ and $H\cap K\leq H_{sG}$, where $H_{sG}$ denotes the subgroup generated by all those subgroups of $H$ which are S-quasinormal in $G$. Let $\mathfrak{F}$ be a formation. A subgroup $H$ of $G$ is called to be $\mathfrak{F}$-quasinormal \cite{Mia} in $G$ if $G$ has a quasinormal subgroup $K$ such that $HK$ is quasinormal in $G$ and $(H\cap K)H_G/H_G\leq Z_{\mathfrak{F}}(G/H_G)$. A subgroup $H$ of $G$ is said to be $\mathfrak{F}_s$-quasinormal \cite{Hua1} in $G$ if $G$ has a normal subgroup $K$ such that $HK$ is S-quasinormal in $G$ and $(H\cap K)H_G/H_G\leq Z_{\mathfrak{F}}(G/H_G)$.\par
\medskip
\noindent\textbf{Lemma 7.2.} Let $H$ be a subgroup of $G$. Then $H$ satisfies partial $\Pi$-property in $G$, if one of the following holds:\par
(1) $H$ is a partial CAP-subgroup of $G$.\par
(2) $H$ is S-embedded in $G$.\par
(3) $H$ is $\mathfrak{U}$-quasinormal in $G$.\par
(4) $H$ is $\mathfrak{U}_s$-quasinormal in $G$.\par
\noindent\textit{Proof.} Statement (1) directly follows from definitions of partial CAP-subgroups and partial $\Pi$-property.\par
(2) Suppose that $H$ is S-embedded in $G$. Then $G$ has a normal subgroup $K$ such that $HK$ is S-quasinormal in $G$ and $H\cap K\leq H_{sG}$. By Lemma 7.1(3), $HK$ satisfies $\Pi$-property, and thus satisfies partial $\Pi$-property in $G$. It follows from Lemma 2.1(3) that $HK/K$ satisfies partial $\Pi$-property in $G/K$. Note that $H_{sG}$ is S-quasinormal in $G$ by \cite[Corollary 1]{Sch}. Hence $H_{sG}$ also satisfies $\Pi$-property in $G$ by Lemma 7.1(3). This implies that for every $G$-chief factor $A/B$ below $K$, $|G:N_G(H_{sG}B\cap A)|$ is a $\pi((H_{sG}B\cap A)/B)$-number. Since $H\cap K=H_{sG}\cap K$, we have that $|G:N_G(HB\cap A)|$ is a $\pi((HB\cap A)/B)$-number. Then by Lemma 2.1(5), $H$ satisfies partial $\Pi$-property in $G$.\par
(3) Suppose that $H$ is $\mathfrak{U}$-quasinormal in $G$. Then $G$ has a quasinormal subgroup $K$ such that $HK$ is quasinormal in $G$ and $(H\cap K)H_G/H_G\leq Z_{\mathfrak{U}}(G/H_G)$. It follows from \cite[Lemma 2.2(2)]{Mia} that $H/H_G$ is $\mathfrak{U}$-quasinormal in $G/H_G$. If $H_G\neq 1$, then by induction, $H/H_G$ satisfies partial $\Pi$-property in $G/H_G$. Hence $H$ satisfies partial $\Pi$-property in $G$ by Lemma 2.1(5). We may, therefore, assume that $H_G=1$. Since $HK$ is quasinormal in $G$, $HK^G$ is also quasinormal in $G$. By Lemma 7.1(3) and Lemma 2.1(3), $HK^G/K^G$ satisfies partial $\Pi$-property in $G/K^G$. As $K$ is quasinormal in $G$, $K^G/K_G\leq Z_\infty(G/K_G)\leq Z_\mathfrak{U}(G/K_G)$ by \cite[Theorem]{Mai}. Then it is easy to see that for every $G$-chief factor $A/B$ with $K_G\leq B\leq A\leq K^G$, $|G:N_G(HB\cap A)|$ is a $\pi((HB\cap A)/B)$-number. By Lemma 2.1(5), $HK_G/K_G$ satisfies partial $\Pi$-property in $G/K_G$. Since $H\cap K$ is $\mathfrak{U}$-hypercentrally embedded in $G$, $H\cap K$ satisfies $\Pi$-property in $G$ by Lemma 7.1(2). Then for every $G$-chief factor $A/B$ below $K_G$, $|G:N_G((H\cap K)B\cap A)|$ is a $\pi(((H\cap K)B\cap A)/B)$-number, and so $|G:N_G(HB\cap A)|$ is a $\pi((HB\cap A)/B)$-number. By Lemma 2.1(5) again, $H$ satisfies partial $\Pi$-property in $G$.\par
Statement (4) can be handled similarly as (3).\par
\medskip
Recall that a subgroup $H$ of $G$ is called to be $\Pi$-normal \cite{Li1} in $G$ if $G$ has a subnormal subgroup $K$ such that $G=HK$ and $H\cap K\leq I\leq H$, where $I$ satisfies $\Pi$-property in $G$. A subgroup $H$ of $G$ is said to be $\mathfrak{U}_c$-normal \cite{Ahm} in $G$ if $G$ has a subnormal subgroup $K$ such that $G=HK$ and $(H\cap K)H_G/H_G\leq Z_{\mathfrak{U}}(G/H_G)$. A subgroup $H$ of $G$ is called to be weakly S-permutable \cite{Ski3} in $G$ if $G$ has a subnormal subgroup $K$ such that $G=HK$ and $H\cap K\leq H_{sG}$, where $H_{sG}$ denotes the subgroup generated by all those subgroups of $H$ which are S-quasinormal in $G$. A subgroup $H$ of $G$ is said to be weakly S-semipermutable \cite{Li4} in $G$ if $G$ has a subnormal subgroup $K$ such that $G=HK$ and $H\cap K\leq H_{ssG}$, where $H_{ssG}$ denotes an S-semipermutable subgroup of $G$ contained in $H$. A subgroup $H$ of $G$ is called to be weakly SS-permutable \cite{He} in $G$ if $G$ has a subnormal subgroup $K$ such that $G=HK$ and $H\cap K\leq H_{ss}$, where $H_{ss}$ denotes an SS-quasinormal subgroup of $G$ contained in $H$. A subgroup $H$ of $G$ is said to be $\tau$-quasinormal \cite{Luk1} in $G$ if $HG_p=G_pH$ for every $G_p\in Syl_p(G)$ such that $(p,|H|)=1$ and $(|H|,|{G_p}^G|)\neq 1$. A subgroup $H$ of $G$ is called to be weakly $\tau$-quasinormal \cite{Luk1} in $G$ if $G$ has a subnormal subgroup $K$ such that $G=HK$ and $H\cap K\leq H_{\tau G}$, where $H_{\tau G}$ denotes the subgroup generated by all those subgroups of $H$ which are $\tau$-quasinormal in $G$.\par
\medskip
\noindent\textbf{Lemma 7.3.} Let $H$ be a $p$-subgroup of $G$. Then $H$ satisfies partial $\Pi$-property in $G$, if one of the following holds:\par
(1) $H$ is $\Pi$-normal in $G$.\par
(2) $H$ is $\mathfrak{U}_c$-normal in $G$.\par
(3) $H$ is weakly S-permutable in $G$.\par
(4) $H$ is weakly S-semipermutable in $G$.\par
(5) $H$ is weakly SS-permutable in $G$.\par
(6) $H$ is weakly $\tau$-quasinormal in $G$.\par
\noindent\textit{Proof.} (1) Suppose that $H$ is $\Pi$-normal in $G$. Then $G$ has a subnormal subgroup $K$ such that $G=HK$ and $H\cap K\leq I\leq H$, where $I$ satisfies $\Pi$-property in $G$. Since $|G:K|$ is a $p$-number, $O^p(G)\leq K$. It follows that for every $G$-chief factor $A/B$ below $O^p(G)$, $|G:N_G(IB\cap A)|$ is a $p$-number. As $H\cap K=I\cap K$, we have that $|G:N_G(HB\cap A)|$ is a $p$-number. Clearly, $HO^p(G)/O^p(G)$ satisfies partial $\Pi$-property in $G/O^p(G)$. Then by Lemma 2.1(5), $H$ satisfies partial $\Pi$-property in $G$.\par
Statements (2)-(5) directly follow from (1) and the fact that a $\mathfrak{U}_c$-normal (resp. weakly S-permutable, weakly S-semipermutable, weakly SS-permutable) subgroup of $G$ is $\Pi$-normal in $G$ by Lemma 7.1.\par
(6) Suppose that $H$ is weakly $\tau$-quasinormal in $G$. Then $G$ has a subnormal subgroup $K$ such that $G=HK$ and $H\cap K\leq H_{\tau G}$. If $O_{p'}(G)\neq 1$, then by \cite[Lemma 2.4(4)]{Luk1}, $HO_{p'}(G)/O_{p'}(G)$ is weakly $\tau$-quasinormal in $G/O_{p'}(G)$. By induction, $HO_{p'}(G)/O_{p'}(G)$ satisfies partial $\Pi$-property in $G/O_{p'}(G)$. Then by Lemma 2.1(5), $H$ satisfies partial $\Pi$-property in $G$. We may, therefore, assume that $O_{p'}(G)=1$. By \cite[Lemma 2.3(1)]{Luk1}, $H_{\tau G}$ is $\tau$-quasinormal in $G$. It follows that $H_{\tau G}G_q=G_qH_{\tau G}$ for every $G_q\in Syl_q(G)$ with $q\in \pi(G)$ such that $q\neq p$ and $p\mid |{G_q}^G|$. As $O_{p'}(G)=1$, we have that $p\mid |{G_q}^G|$ for every $q\neq p$. Hence $H_{\tau G}$ is S-semipermutable in $G$, and so $H$ is weakly S-semipermutable in $G$. By (4), $H$ satisfies partial $\Pi$-property in $G$.\par
\medskip
Now our attention is restricted to the solvable universe. Recall that a subgroup $H$ of $G$ is said to be S-quasinormally embedded \cite{Bal7} in $G$ if every Sylow subgroup of $H$ is a Sylow subgroup of some S-quasinormal subgroup of $G$. A subgroup $H$ of $G$ is called to be S-conditionally permutable \cite{Hua} in $G$ if $H$ permutes with at least one Sylow $p$-subgroup of $G$ for every $p\in\pi(G)$. A subgroup $H$ of $G$ is said to be S-C-permutably embedded \cite{Che1} in $G$ if every Sylow subgroup of $H$ is a Sylow subgroup of some S-conditionally permutable subgroup of $G$.\par
\medskip
\noindent\textbf{Lemma 7.4.} Let $H$ be a subgroup of $G$ contained in a solvable normal subgroup $N$ of $G$. Then $H$ satisfies $\Pi$-property, and thus satisfies partial $\Pi$-property in $G$, if one of the following holds:\par
(1) $H$ is S-quasinormally embedded in $G$.\par
(2) $H$ is S-conditionally permutable in $G$.\par
(3) $H$ is S-C-permutably embedded in $G$.\par
\noindent\textit{Proof.} Statements (1)-(2) were proved in \cite[Proposition 2.5]{Li1}.\par
(3) Suppose that $H$ is S-C-permutably embedded in $G$. Let $L/K$ be a $G$-chief factor. Then we only need to prove that $|G:N_G(HK\cap L)|$ is a $\pi((HK\cap L)/K)$-number. By \cite[Lemma 2.2(1)]{Che1}, we may assume that $K=1$ by induction. If $N\cap L=1$, then $H\cap L=1$, there is nothing to prove. Now suppose that $L\leq N$. Then $L$ is a $p$-group with $p\in \pi(G)$. Since $H$ is S-C-permutably embedded in $G$, $G$ has an S-conditionally permutable subgroup $X$ such that some Sylow $p$-subgroup of $H$ is a Sylow $p$-subgroup of $X$. This implies that $H\cap L=X\cap L$. As $X$ is S-conditionally permutable in $G$, for every $q\in \pi(G)$ with $p\neq q$, $G$ has a Sylow $q$-subgroup $G_q$ such that $X$ permutes with $G_q$. It follows that $H\cap L=X\cap L=XG_q\cap L\unlhd XG_q$, and thereby $G_q\leq N_G(H\cap L)$. Hence $|G:N_G(H\cap L)|$ is a $p$-number. This shows that $H$ satisfies $\Pi$-property in $G$.\par
\medskip
Recall that a subgroup $H$ of $G$ is said to be weakly S-embedded \cite{Li3} in $G$ if $G$ has a normal subgroup $K$ such that $HK$ is S-quasinormal in $G$ and $H\cap K\leq H_{seG}$, where $H_{seG}$ denotes the subgroup generated by all those subgroups of $H$ which are S-quasinormally embedded in $G$. A subgroup $H$ of $G$ is called to be weakly S-permutably embedded \cite{Li5} in $G$ if $G$ has a subnormal subgroup $K$ such that $G=HK$ and $H\cap K\leq H_{se}$, where $H_{se}$ denotes an S-quasinormally embedded subgroup of $G$ contained in $H$.\par
\medskip
\noindent\textbf{Lemma 7.5.} Let $H$ be a $p$-subgroup of $G$ contained in a solvable normal subgroup $N$ of $G$. Then $H$ satisfies partial $\Pi$-property in $G$, if one of the following holds:\par
(1) $H$ is weakly S-embedded in $G$.\par
(2) $H$ is weakly S-permutably embedded in $G$.\par
\noindent\textit{Proof.} (1) Suppose that $H$ is weakly S-embedded in $G$. Then $G$ has a normal subgroup $K$ such that $HK$ is S-quasinormal in $G$ and $H\cap K\leq H_{seG}$. By \cite[Lemma 2.4(1)]{Li3}, $H/H_G$ is weakly S-embedded in $G$. If $H_G\neq 1$, then by induction, $H/H_G$ satisfies partial $\Pi$-property in $G/H_G$. Hence $H$ satisfies partial $\Pi$-property in $G$ by Lemma 2.1(5). We may, therefore, assume that $H_G=1$. By \cite[Lemma 2.4(3)]{Li3}, $HO_{p'}(G)/O_{p'}(G)$ is weakly S-embedded in $G/O_{p'}(G)$. Similarly as above, we may assume that $O_{p'}(G)=1$. Let $H_1,H_2,\cdots ,H_n$ be all subgroups of $H$ which are S-quasinormally embedded in $G$. Then there exist S-quasinormal subgroups $X_1,X_2,\cdots ,X_n$ of $G$ with $H_i\in Syl_p(X_i)$ for $1\leq i\leq n$. Suppose that ${(X_k)}_G\neq 1$ for some integer $k$. Let $N_k$ be a minimal normal subgroup of $G$ contained in $(X_k)_G$. If $N\cap N_k=1$, then $H_k\cap N_k=1$, and so $p\nmid |N_k|$. Hence $N_k\leq O_{p'}(G)=1$, a contradiction. Thus $N_k\leq N$. Since $O_{p'}(G)=1$, $N_k$ is a $p$-group. This implies that $N_k\leq (H_k)_G=1$, which is impossible. Consequently, we get that ${(X_i)}_G=1$ for $1\leq i\leq n$. By \cite[Proposition A]{Sch}, $X_i$ is nilpotent  for $1\leq i\leq n$. It follows from \cite[Proposition B]{Sch} that $H_i$ is S-quasinormal in $G$ for $1\leq i\leq n$. Therefore, by \cite[Corollary 1]{Sch}, $H_{seG}=\langle H_1,H_2,\cdots ,H_n\rangle$ is S-quasinormal in $G$. This shows that $H$ is S-embedded in $G$. Hence $H$ satisfies partial $\Pi$-property in $G$ by Lemma 7.2(2).\par
Statement (2) directly follows from Lemma 7.3(1) and the fact that a weakly S-permutably embedded subgroup of $G$ is $\Pi$-normal in $G$ by Lemma 7.4(1).\par
\medskip
By the above lemmas, one can see that a lot of previous results can be deduced from our theorems. Interested readers may refer to the relevant literature for further details.\par
\bibliographystyle{plain}
\bibliography{expbib}
\end{document}